\documentclass{article}%
\usepackage{amsmath}
\usepackage{amsfonts}
\usepackage{amssymb}
\usepackage{graphicx}%
\newcommand{\bpartial}{\mathop{\partial\kern -4pt\raisebox{.8pt}{$|$}}}
\newcommand{\bra}{\mathopen{[\kern-1.6pt[}}
\newcommand{\ket}{\mathclose{]\kern-1.5pt]}}
\newcommand{\bbra}{\mathopen{[\kern-2.2pt[\kern-2.3pt[}}
\newcommand{\bket}{\mathclose{]\kern-2.1pt]\kern-2.3pt]}}
\newcommand{\slg}{\mbox{\bfseries\slshape g}}

\newcommand{\slx}{\mbox{\bfseries\slshape x}}

\makeindex
\begin{document}

\title{Clifford and Extensor Calculus and the Riemann and Ricci Extensor Fields of
Deformed Structures $(M,\nabla^{\prime},\mathbf{\eta})$ and $(M,\nabla,%
\slg
)$\thanks{I\textit{nt. J. Geom Meth. Mod. Phys.} \textbf{4} (7), 1159-1172
(2007)}}
\author{{\footnotesize V. V. Fern\'{a}ndez}$^{1}${\footnotesize , W. A. Rodrigues
Jr}.$^{1}$\ {\footnotesize A. M. Moya}$^{2}${\footnotesize  and  R. da Rocha
Jr.}$^{3}$ \\$^{1}\hspace{-0.1cm}${\footnotesize Institute of Mathematics, Statistics and
Scientific Computation}\\{\footnotesize \ IMECC-UNICAMP CP 6065}\\{\footnotesize \ 13083-859 Campinas, SP, Brazil }\\{\footnotesize e-mail:walrod@ime.unicamp.br virginvelfe@accessplus.com.ar}\\{\footnotesize \ }$^{2}${\footnotesize Department of Mathematics, University
of Antofagasta, Antofagasta, Chile} \\{\footnotesize e-mail: mmoya@uantof.cl}\\$^{3}${\footnotesize Instituto de F\'{\i}sica Te\'{o}rica}\\{\footnotesize Universidade Estadual Paulista}\\{\footnotesize Rua Pamplona 145, 01405-900 S\~{a}o Paulo, SP, Brazil}\\{\footnotesize and}\\{\footnotesize Institute of Physics G. Wataghin, UNICAMP}\\{\footnotesize e-mail: roldao@ifi.unicamp.br}}
\maketitle

\begin{abstract}
Here (the last paper in a series of four) we end our presentation of the
basics of a systematical approach to the differential geometry of a smooth
manifold $M$ (supporting a metric field $%
\slg
$ \ and a general connection $\nabla$) which uses the geometric algebras of
multivector and extensors (fields) developed in previous papers. The theory of
the Riemann and Ricci fields of the the triple $(M,\nabla,%
\slg
)$ is investigated to for each particular open set $U\subset M$ through the
introduction of a geometric structure on $U$, i.e., a triple $(U,\gamma,g)$
where $\gamma$ is a general connection field on $U$ and $g$ is a metric
extensor field associated to $%
\slg
$. The relation between geometrical structures related by gauge extensor
fields is clarified. These geometries may be said to be deformations one of
each other. Moreover we study the important case of a class of deformed
Levi-Civita geometrical structures and prove key theorems about them that are
important in the formulation of geometric theories of the gravitational field.

\end{abstract}

\section{Introduction}

This is the last paper in a series of four where we systematically applied the
geometric algebras of multivector and extensors \cite{1} to the differential
geometry of an arbitrary smooth manifold $M$ equipped with a metric field $%
\slg
$ and an arbitrary connection $\nabla$ \cite{4,5}. Here we study in Section 2
the theory of the Riemann and Ricci fields of a triple $(M,\nabla,%
\slg
)$ on an arbitrary open set $U\subset M$ through the introduction of a
geometric structure on $U$, i.e., a triple $(U,\gamma,g)$ where $\gamma$ is a
general connection field on $U$ and $g$ is a metric extensor field associated
to $%
\slg
$. We concentrate our study on the relation between geometrical structures
$(U,\gamma^{\prime},\eta)$ and $(U,\gamma,g)$ related by gauge extensor
fields. These geometries may be said to be \textit{deformations} one of each
other. In Section 3 we study in details the important case of a class of
deformed Levi-Civita geometrical structures and prove some key theorems about
these structures that are important in the formulation of geometric theories
of the gravitational field. We show that if we start with a flat geometry for
$(M,\nabla^{\prime},\mathbf{\eta})$ with \textit{null curvature} and
\textit{null torsion} tensors, represented in $U$ by the geometric structure
$(U,\gamma^{\prime},\eta)$, the deformed geometry $(U,\gamma,g)$, where the
covariant derivative of $\gamma$ is as defined in \cite{5}, can only be
associated with a triple $(M,\nabla,%
\slg
)$ which has also \textit{null} curvature tensor but in general has a
\textit{non null} torsion tensor. Non null curvature for $(U,\gamma,g)$ only
occurs if $(U,\gamma^{\prime},\eta)$ has non null curvature. We introduce also
the concept of deformed covariant derivatives associated to a diffeomorphism
$\mathtt{h}:M\rightarrow M$ (\texttt{h}$(U)\subset U$) and show that there is
a deformation tensor $h$ associate to such a diffeomorphism such that if
$(U,\gamma^{\prime},\eta)$, has null curvature and torsion tensors, the
deformed geometry $(U,\gamma=\mathtt{h}^{\ast}\gamma^{\prime},%
\slg
=\mathtt{h}^{\ast}%
\mbox{\boldmath{$\eta$}}%
)$ also has null curvature and torsion tensors\footnote{Applications of these
ideas in geometrical theories of the gravitational field has been presented in
\cite{notterod}.}. In section 4 we present our conclusions.

\section{Geometric Structure}

Let $U$ be an open subset of $U_{o},$ and let $(U,\gamma,g)$ be a
\emph{geometric structure} on $U.$ We recall, from a previous paper in the
series \cite{5}, that this means that the \emph{connection field} $\gamma$ is
\emph{compatible} with the \emph{metric field} $g,$ i.e., $\gamma
_{a+(g)}=\dfrac{1}{2}g^{-1}\circ(a\cdot\partial_{o}g).$ From a logically
equivalent point of view, this implies that the pair of $a$-\emph{DCDO's}
associated to the \emph{parallelism structure} $(U,\gamma,),$ namely
$(\mathcal{D}_{a}^{+},\mathcal{D}_{a}^{-}),$ is \emph{compatible }with $g,$
i.e., $\mathcal{D}_{a}^{++}g=0$ (or equivalently, $\mathcal{D}_{a}^{--}%
g^{-1}=0$).

The $g$-compatibility of $(U,\gamma,g)$ implies that there exists an unique
smooth $(1,2)$-extensor field, namely $\omega,$ such that $\mathcal{D}_{a}%
^{+}$ is given by
\begin{equation}
\mathcal{D}_{a}^{+}X=a\cdot\partial_{o}X+\frac{1}{2}\underline{g}^{-1}%
\circ(a\cdot\partial_{o}\underline{g})(X)+\omega(a)\underset{g}{\times}X.
\label{GS.1}%
\end{equation}
Also, the relationship between $\mathcal{D}_{a}^{+}$ and $\mathcal{D}_{a}^{-}$
is given by
\begin{equation}
\underline{g}(\mathcal{D}_{a}^{+}X)=\mathcal{D}_{a}^{-}\underline{g}(X).
\label{GS.2}%
\end{equation}

\subsection{Riemann and Ricci Fields}

Associated to $(U,\gamma,g),$ the smooth vector elementary $3$-extensor field
on $U,$ namely $\rho,$ defined by
\begin{equation}
\rho(a,b,c)=[\mathcal{D}_{a}^{+},\mathcal{D}_{b}^{+}]c-\mathcal{D}_{[a,b]}%
^{+}c,\text{ for all }a,b,c\in\mathcal{V}(U) \label{RRF.1}%
\end{equation}
is said to be the \emph{curvature field }of $(U,\gamma,g).$

We present here the basic properties satisfied by the curvature field.\vspace
{0.1in}

\textbf{i. }$\rho$ has a skew-symmetry property by interchanging the first and
second vector variables, i.e.,
\begin{equation}
\rho(a,b,c)=-\rho(b,a,c). \label{RRF.2}%
\end{equation}
It is an obvious result completely general\footnote{For any parallelism
structure $(U,\gamma)$ whose pair of $a$-\emph{DCDO's }associated is
$(\nabla_{a}^{+},\nabla_{a}^{-})$, the curvature field $\rho$ is defined by
$\rho(a,b,c)=[\nabla_{a}^{+},\nabla_{b}^{+}]c-\nabla_{\lbrack a,b]}^{+}c.$
Then, because of the skew-symmetry properties of $[\nabla_{a}^{+},\nabla
_{b}^{+}]$ and $[a,b],$ it follows that $\rho(a,b,c)=-\rho(b,a,c).$}.

\textbf{ii.} For all $a,b,c\in\mathcal{V}(U)$ it holds
\begin{equation}
\rho(a,b,c)\underset{g}{\cdot}c=0. \label{RRF.2a}%
\end{equation}

The proof is as follows. We use twice the Ricci-like theorem for
$\mathcal{D}_{a}^{+}$ we can write
\begin{align}
a\cdot\partial_{o}b\cdot\partial_{o}(c\underset{g}{\cdot}c)  &  =a\cdot
\partial_{o}((\mathcal{D}_{b}^{+}c)\underset{g}{\cdot}c+c\underset{g}{\cdot
}(\mathcal{D}_{b}^{+}c))=2a\cdot\partial_{o}((\mathcal{D}_{b}^{+}%
c)\underset{g}{\cdot}c)\nonumber\\
&  =2\mathcal{D}_{a}^{+}(\mathcal{D}_{b}^{+}c)\underset{g}{\cdot
}c+2(\mathcal{D}_{b}^{+}c)\underset{g}{\cdot}(\mathcal{D}_{a}^{+}c),
\label{RRF.N1}%
\end{align}
and, by interchanging the letters $a$ and $b,$ we get
\begin{equation}
b\cdot\partial_{o}a\cdot\partial_{o}(c\underset{g}{\cdot}c)=2\mathcal{D}%
_{b}^{+}(\mathcal{D}_{a}^{+}c)\underset{g}{\cdot}c+2(\mathcal{D}_{a}%
^{+}c)\underset{g}{\cdot}(\mathcal{D}_{b}^{+}c). \label{RRF.N2}%
\end{equation}
Now, subtracting Eq.(\ref{RRF.N2}) from Eq.(\ref{RRF.N1}), we have
\[
\lbrack a\cdot\partial_{o},b\cdot\partial_{o}](c\underset{g}{\cdot
}c)=2[\mathcal{D}_{a}^{+},\mathcal{D}_{b}^{+}]c\underset{g}{\cdot}c.
\]
But, by recalling the identity $[a\cdot\partial_{o},b\cdot\partial
_{o}]X=[a,b]\cdot\partial_{o}X,$ we can write
\[
\lbrack a,b]\cdot\partial_{o}(c\underset{g}{\cdot}c)=2[\mathcal{D}_{a}%
^{+},\mathcal{D}_{b}^{+}]c\underset{g}{\cdot}c,
\]
and, by using once more the Ricci-like theorem for $\mathcal{D}_{a}^{+},$ we
get
\begin{equation}
2(\mathcal{D}_{[a,b]}^{+}c)\underset{g}{\cdot}c=2[\mathcal{D}_{a}%
^{+},\mathcal{D}_{b}^{+}]c\underset{g}{\cdot}c. \label{RRF.N3}%
\end{equation}
The required result immediately follows from Eq.(\ref{RRF.N3}).

\textbf{iii. }For all $a,b,c,d\in\mathcal{V}(U)$ it holds
\begin{equation}
\rho(a,b,c)\underset{g}{\cdot}d=-\rho(a,b,d)\underset{g}{\cdot}c.
\label{RRF.2b}%
\end{equation}

To show Eq.\ref{RRF.2b} we use Eq.(\ref{RRF.2a}) and write
\begin{align*}
\rho(a,b,c+d)\underset{g}{\cdot}(c+d)  &  =0\\
\rho(a,b,c)\underset{g}{\cdot}c+\rho(a,b,c)\underset{g}{\cdot}d+\rho
(a,b,d)\underset{g}{\cdot}c+\rho(a,b,d)\underset{g}{\cdot}d  &  =0\\
\rho(a,b,c)\underset{g}{\cdot}d+\rho(a,b,d)\underset{g}{\cdot}c  &  =0,
\end{align*}
and the expected result follows.

We emphasize that Eq.(\ref{RRF.2a}) and Eq.(\ref{RRF.2b}) are \emph{logically
equivalent} to each other.

\textbf{Proposition 1}. There exists a smooth $(2,2)$-extensor field on $U,$
namely $B\mapsto R(B),$ such that for all $a,b,c,d\in\mathcal{V}(U)$
\begin{equation}
-\rho(a,b,c)\underset{g}{\cdot}d=R(a\wedge b)\cdot(c\wedge d). \label{RRF.3}%
\end{equation}

Such $B\mapsto R(B)$ is given by
\begin{equation}
R(B)=-\frac{1}{4}B\cdot(\partial_{a}\wedge\partial_{b})\partial_{c}%
\wedge\partial_{d}\rho(a,b,c)\underset{g}{\cdot}d, \label{RRF.3a}%
\end{equation}
i.e.,
\begin{equation}
R(a\wedge b)=-\frac{1}{2}\partial_{c}\wedge\partial_{d}\rho(a,b,c)\underset
{g}{\cdot}d. \label{RRF.3b}%
\end{equation}

\textbf{Proof}

We emphasize that Eq.(\ref{RRF.3a}) and Eq.(\ref{RRF.3b}) are in fact
\emph{logically equivalent} to each other.

Eq.(\ref{RRF.3a}) implies Eq.(\ref{RRF.3b}), i.e.,
\begin{align*}
R(a\wedge b)  &  =-\frac{1}{4}(a\wedge b)\cdot(\partial_{p}\wedge\partial
_{q})\partial_{c}\wedge\partial_{d}\rho(p,q,c)\underset{g}{\cdot}d\\
&  =-\frac{1}{4}\det\left[
\begin{array}
[c]{cc}%
a\cdot\partial_{p} & a\cdot\partial_{q}\\
b\cdot\partial_{p} & b\cdot\partial_{q}%
\end{array}
\right]  \partial_{c}\wedge\partial_{d}\rho(p,q,c)\underset{g}{\cdot}d\\
&  =-\frac{1}{4}(a\cdot\partial_{p}b\cdot\partial_{q}-a\cdot\partial_{q}%
b\cdot\partial_{p})\partial_{c}\wedge\partial_{d}\rho(p,q,c)\underset{g}%
{\cdot}d\\
&  =-\frac{1}{2}(a\cdot\partial_{p}b\cdot\partial_{q})\partial_{c}%
\wedge\partial_{d}\rho(p,q,c)\underset{g}{\cdot}d\\
&  =-\frac{1}{2}\partial_{c}\wedge\partial_{d}\rho(a,b,c)\underset{g}{\cdot}d.
\end{align*}
We have used the skew-symmetry property given by Eq.(\ref{RRF.2}).

We now show that Eq.(\ref{RRF.3b}) implies Eq.(\ref{RRF.3a}). We can write
\begin{align*}
R(B)  &  =R(\frac{1}{2}B\cdot(\partial_{a}\wedge\partial_{b})a\wedge
b)=\frac{1}{2}B\cdot(\partial_{a}\wedge\partial_{b})R(a\wedge b)\\
&  =-\frac{1}{4}B\cdot(\partial_{a}\wedge\partial_{b})\partial_{c}%
\wedge\partial_{d}\rho(a,b,c)\underset{g}{\cdot}d.
\end{align*}

In order to prove that Eq.(\ref{RRF.3}) can be deduced by using
Eq.(\ref{RRF.3b}), we have
\begin{align*}
R(a\wedge b)\cdot(c\wedge d)  &  =-\frac{1}{2}(c\wedge d)\cdot(\partial
_{p}\wedge\partial_{q})\rho(a,b,p)\underset{g}{\cdot}q\\
&  =-\frac{1}{2}\det\left[
\begin{array}
[c]{cc}%
c\cdot\partial_{p} & c\cdot\partial_{q}\\
d\cdot\partial_{p} & d\cdot\partial_{q}%
\end{array}
\right]  \rho(a,b,p)\underset{g}{\cdot}q\\
&  =-\frac{1}{2}(c\cdot\partial_{p}d\cdot\partial_{q}-c\cdot\partial_{q}%
d\cdot\partial_{p})\rho(a,b,p)\underset{g}{\cdot}q\\
&  =-c\cdot\partial_{p}d\cdot\partial_{q}\rho(a,b,p)\underset{g}{\cdot}%
q=-\rho(a,b,c)\underset{g}{\cdot}d.
\end{align*}
We have used the skew-symmetry property given by Eq.(\ref{RRF.2b}%
).$\blacksquare$

Such $B\mapsto R(B)$ will be called the \emph{Riemann field} for
$(U,\gamma,g).$

\textbf{Proposition 2.}There exists a smooth $(1,1)$-extensor field on $U,$
namely $b\mapsto R(b),$ such that for all $b,c\in\mathcal{V}(U)$
\begin{equation}
R(b)\cdot c=\partial_{a}\cdot\rho(a,b,c). \label{RRF.4}%
\end{equation}

Such $b\mapsto R(b)$ is given by
\begin{equation}
R(b)=\partial_{c}\partial_{a}\cdot\rho(a,b,c). \label{RRF.4a}%
\end{equation}

\textbf{Proof}

We have indeed that
\[
R(b)\cdot c=c\cdot\partial_{p}\partial_{a}\cdot\rho(a,b,p)=\partial_{a}%
\cdot\rho(a,b,c).\blacksquare
\]

Such $b\mapsto R(b)$ will be called the \emph{Ricci field} for $(U,\gamma,g).$
The smooth scalar field on $U,$ namely $R,$ defined by
\begin{equation}
R=g^{-1}(\partial_{b})\cdot R(b), \label{RRF.5}%
\end{equation}
is called the \emph{Ricci curvature scalar field} for $(U,\gamma,g).$

We present now some noticeable properties involving the Riemann and the Ricci fields.

Let us take $a,b,c,d\in\mathcal{V}(U).$ By using Eq.(\ref{RRF.3}) and the
identity $X\cdot(Y\wedge Z)=(X\llcorner\widetilde{Z})\cdot Y,$ where
$X,Y,Z\in\bigwedge\mathcal{U}_{o},$ we can write
\[
g\circ\rho(a,b,c)\cdot d=-R(a\wedge b)\cdot(c\wedge d)=R(a\wedge
b)\cdot(d\wedge c)=(R(a\wedge b)\llcorner c)\cdot d,
\]
hence, by the non-degeneracy of the scalar product, we get
\begin{equation}
g\circ\rho(a,b,c)=R(a\wedge b)\times c. \label{RRF.6a}%
\end{equation}

By using Eq.(\ref{RRF.4}), Eq.(\ref{RRF.6a}) and the identities $X\cdot
(Y\llcorner Z)=(X\wedge\widetilde{Z})\cdot Y$ and $(X\wedge Y)\cdot
Z=Y\cdot(\widetilde{X}\lrcorner Z),$ where $X,Y,Z\in\bigwedge\mathcal{U}_{o},$
we have
\begin{align*}
R(b)\cdot c  &  =\partial_{a}\cdot g^{-1}(R(a\wedge b)\llcorner c)=g^{-1}%
(\partial_{a})\cdot R(a\wedge b)\llcorner c\\
&  =(g^{-1}(\partial_{a})\wedge c)\cdot R(a\wedge b)=c\cdot(g^{-1}%
(\partial_{a})\lrcorner R(a\wedge b)),
\end{align*}
hence, by the non-degeneracy of the scalar product, we get
\begin{equation}
R(b)=g^{-1}(\partial_{a})\lrcorner R(a\wedge b). \label{RRF.6b}%
\end{equation}

A straightforward calculation by using Eq.(\ref{RRF.5}) and Eq.(\ref{RRF.6b})
allows us to get
\begin{align}
R  &  =g^{-1}(\partial_{b})\cdot(g^{-1}(\partial_{a})\lrcorner R(a\wedge
b))\nonumber\\
&  =(g^{-1}(\partial_{a})\wedge g^{-1}(\partial_{b}))\cdot R(a\wedge
b),\nonumber\\
R  &  =\underline{g}^{-1}(\partial_{a}\wedge\partial_{b})\cdot R(a\wedge b).
\label{RRF.6c}%
\end{align}

\textbf{Proposition 3. }For all $a,b\in\mathcal{V}(U)$ and $X\in
\mathcal{M}(U)$ the following property holds
\begin{equation}
\lbrack\mathcal{D}_{a}^{+},\mathcal{D}_{b}^{+}]X=\mathcal{D}_{[a,b]}%
^{+}X+\underline{g}^{-1}(R(a\wedge b)\underset{g^{-1}}{\times}\underline
{g}(X)). \label{RRF.7}%
\end{equation}

\bigskip\textbf{Proof}

In order to prove this noticeable property we need to check it only for smooth
scalar fields, and for the scalar multiplication of smooth simple $k$-vector
fields by smooth scalar fields, i.e.,
\begin{align}
\lbrack\mathcal{D}_{a}^{+},\mathcal{D}_{b}^{+}]f  &  =\mathcal{D}_{[a,b]}%
^{+}f+\underline{g}^{-1}(R(a\wedge b)\underset{g^{-1}}{\times}\underline
{g}(f)),\label{RRF.7a}\\
\lbrack\mathcal{D}_{a}^{+},\mathcal{D}_{b}^{+}](fv_{1}\wedge\ldots\wedge
v_{k})  &  =\mathcal{D}_{[a,b]}^{+}(fv_{1}\wedge\ldots\wedge v_{k})\nonumber\\
&  +\underline{g}^{-1}(R(a\wedge b)\underset{g^{-1}}{\times}\underline
{g}(fv_{1}\wedge\ldots\wedge v_{k})). \label{RRF.7b}%
\end{align}

The proof of the first statement is trivial. To prove the second statement we
should use the Leibniz rules for scalar multiplication and exterior product,
and Eq.(\ref{RRF.6a}).$\blacksquare$

\subsection{Gauge Extensor Fields}

Let $h$ be a gauge metric field for $g$. As we know from the previous paper in
this series, this means that there exists a smooth $(1,1)$-extensor field,
namely $h$, such that $g=h^{\dagger}\circ\eta\circ h$ where $\eta$ is an
orthogonal metric field with the same signature as $g$.

To the $g$-compatible pair of $a$-\emph{DCDO's }associated to $(U,\gamma,g),$
namely $(_{g}\mathcal{D}_{a}^{+},_{g}\mathcal{D}_{a}^{-}),$ corresponds an
unique $\eta$-compatible pair of $a$-\emph{DCDO's,} namely $(_{\eta
}\mathcal{D}_{a}^{+},_{\eta}\mathcal{D}_{a}^{-}),$ such that the latter is
just a $h$-\emph{deformation }of the former, i.e.,
\begin{align}
_{\eta}\mathcal{D}_{a}^{+}X  &  =\underline{h}(_{g}\mathcal{D}_{a}%
^{+}\underline{h}^{-1}(X)),\label{GEF.1a}\\
_{\eta}\mathcal{D}_{a}^{-}X  &  =\underline{h}^{\ast}(_{g}\mathcal{D}_{a}%
^{-}\underline{h}^{\dagger}(X)). \label{GEF.1b}%
\end{align}

The $\eta$-compatibility of $(_{\eta}\mathcal{D}_{a}^{+},_{\eta}%
\mathcal{D}_{a}^{-})$ implies that there must be an unique smooth
$(1,2)$-extensor field, say $\Omega,$ such that $_{\eta}\mathcal{D}_{a}^{+}$
is given by
\begin{equation}
_{\eta}\mathcal{D}_{a}^{+}X=a\cdot\partial_{o}X+\Omega(a)\underset{\eta
}{\times}X. \label{GEF.2a}%
\end{equation}
And, $_{\eta}\mathcal{D}_{a}^{-}$ is related to $_{\eta}\mathcal{D}_{a}^{+}$
by
\begin{equation}
_{\eta}\mathcal{D}_{a}^{-}X=\underline{\eta}(_{\eta}\mathcal{D}_{a}%
^{+}\underline{\eta}(X)). \label{GEF.2b}%
\end{equation}

In agreement with Eq.(\ref{RRF.1}), the curvature field of the $\eta
$-\emph{geometric structure} $(U_{o},(a,b)\mapsto\Omega(a)\underset{\eta
}{\times}b,\eta),$ namely $_{\eta}\rho,$ is given by
\begin{equation}
_{\eta}\rho(a,b,c)=[_{\eta}\mathcal{D}_{a}^{+},_{\eta}\mathcal{D}_{a}%
^{-}]c-\text{ }_{\eta}\mathcal{D}_{[a,b]}^{+}c. \label{GEF.3}%
\end{equation}
But, a straightforward calculation, by using Eq.(\ref{GEF.2a}) and the
remarkable identity $B\underset{\eta}{\times}(X\underset{\eta}{\times
}Y)=(B\underset{\eta}{\times}X)\underset{\eta}{\times}Y+X\underset{\eta
}{\times}(B\underset{\eta}{\times}Y),$ where $B\in\bigwedge^{2}\mathcal{U}%
_{o}$ and $X,Y\in\bigwedge\mathcal{U}_{o},$ yields
\begin{align}
_{\eta}\rho(a,b,c)  &  =(a\cdot\partial_{o}\Omega(b)-b\cdot\partial_{o}%
\Omega(a)-\Omega([a,b])+\Omega(a)\underset{\eta}{\times}\Omega(b))\underset
{\eta}{\times}c,\nonumber\\
&  =((a\cdot\partial_{o}\Omega)(b)-(b\cdot\partial_{o}\Omega)(a)+\Omega
(a)\underset{\eta}{\times}\Omega(b))\underset{\eta}{\times}c. \label{GEF.4}%
\end{align}

By recalling Eq.(\ref{RRF.6a}), the Riemann field for $(U,(a,b)\mapsto
\Omega(a)\underset{\eta}{\times}b,\eta),$ namely $B\mapsto_{\eta}R(B),$ must
satisfy
\begin{equation}
_{\eta}R(a\wedge b)\llcorner c=\eta\circ\text{ }_{\eta}\rho(a,b,c),
\label{GEF.5}%
\end{equation}
by using Eq.(\ref{GEF.4}) and the identity $\underline{\eta}(X\llcorner
\underline{\eta}(Y))=\underline{\eta}(X)\llcorner Y,$ we have
\[
_{\eta}R(a\wedge b)\llcorner c=\underline{\eta}((a\cdot\partial_{o}%
\Omega)(b)-(b\cdot\partial_{o}\Omega)(a)+\Omega(a)\underset{\eta}{\times
}\Omega(b))\llcorner c,
\]
hence,
\begin{equation}
_{\eta}R(a\wedge b)=\underline{\eta}((a\cdot\partial_{o}\Omega)(b)-(b\cdot
\partial_{o}\Omega)(a)+\Omega(a)\underset{\eta}{\times}\Omega(b)).
\label{GEF.6}%
\end{equation}

Next, we relate the curvature field and the Riemann field for $(U,\gamma,g),$
namely $_{g}\rho$ and $B\mapsto$ $_{g}R(B),$ with $_{\eta}\rho$ and $B\mapsto$
$_{\eta}R(B),$ respectively.

In agreement to Eq.(\ref{RRF.1}), and using Eq.(\ref{GEF.1a}), we have
\begin{align}
\text{ }_{g}\rho(a,b,h^{-1}c)  &  =[_{g}\mathcal{D}_{a}^{+},_{g}%
\mathcal{D}_{b}^{+}]h^{-1}c-\text{ }_{g}\mathcal{D}_{[a,b]}^{+}h^{-1}%
c\nonumber\\
&  =h^{-1}([_{\eta}\mathcal{D}_{h^{-1}a}^{+},_{\eta}\mathcal{D}_{h^{-1}b}%
^{+}]c-\text{ }_{\eta}\mathcal{D}_{[a,b]}^{+}c)\nonumber\\
\text{ }  &  =h^{-1}\circ\text{ }_{\eta}\rho(a,b,c) \label{GEF.7}%
\end{align}

Once again Eq.(\ref{RRF.6a}) allows us to write
\[
_{g}R(a\wedge b)\llcorner h^{-1}c=g\circ\text{ }_{g}\rho(a,b,h^{-1}c).
\]
Then, by using the\emph{\ master formula} $g=h^{\dagger}\circ\eta\circ h$ and
Eq.(\ref{GEF.7}), we get
\[
_{g}R(a\wedge b)\llcorner h^{-1}c=h^{\dagger}\circ\eta\circ\text{ }_{\eta}%
\rho(a,b,c)
\]
and by recalling Eq.(\ref{GEF.5}) and using the identity $\underline
{h}^{\dagger}(X\llcorner\underline{h}(Y))=\underline{h}^{\dagger}(X)\llcorner
Y,$ we have
\begin{align}
_{g}R(a\wedge b)\llcorner h^{-1}c  &  =h^{\dagger}(\text{ }_{\eta}R(a\wedge
b)\llcorner c)\nonumber\\
&  =(\underline{h}^{\dagger}\circ\text{ }_{\eta}R(a\wedge b))\llcorner h^{-1}c
\label{deformation 2}%
\end{align}
Then, it immediately follows that
\begin{equation}
_{g}R(B)=\underline{h}^{\dagger}\circ\text{ }_{\eta}R(\text{ }B) \label{GEF.8}%
\end{equation}

Eq.(\ref{GEF.8}) shows that if the geometry $(M,\eta,_{\eta}\mathcal{D})$ is
flat, i.e., $_{\eta}R($ $B)=0$, the geometry $(M,g,_{g}\mathcal{D})$ is also
flat, i.e., $_{g}R(B)=0$.

We introduce now the \emph{gauge Riemann field} and \emph{gauge Ricci field},
namely $B\mapsto\mathcal{R}(B)$ and $b\mapsto\mathcal{R}(b),$ which are
defined as follows
\begin{align}
\mathcal{R}(B)  &  =\underline{\eta\circ h}^{\ast}\circ\text{ }_{g}%
R(B),\label{GEF.9a}\\
\mathcal{R}(b)  &  =\underline{\eta\circ h}^{\ast}\circ\text{ }_{g}R(b).
\label{GEF.9b}%
\end{align}
These smooth extensor fields on $U$ have some interesting properties which are
ready to be used in geometric theories of gravitation.

As we can see, by using Eq.(\ref{GEF.8}) and Eq.(\ref{GEF.6}) into
Eq.(\ref{GEF.9a}), the gauge Riemann field is given by
\begin{equation}
\mathcal{R}(a\wedge b)=(a\cdot\partial_{o}\Omega)(b)-(b\cdot\partial_{o}%
\Omega)(a)+\Omega(a)\underset{\eta}{\times}\Omega(b) \label{GEF.10}%
\end{equation}

By using the identity $X\lrcorner\underline{t}(Y)=\underline{t}(\underline
{t}^{\dagger}(X)\lrcorner Y),$ the \emph{master formula} $g^{-1}=h^{-1}%
\circ\eta\circ h^{\ast},$ and Eq.(\ref{RRF.6b}), a straightforward calculation
yields
\begin{align}
h^{\ast}(\partial_{a})\lrcorner\mathcal{R}(a\wedge b)  &  =h^{\ast}%
(\partial_{a})\lrcorner\underline{\eta\circ h}^{\ast}\circ\text{ }%
_{g}R(a\wedge b)\nonumber\\
&  =\underline{\eta\circ h}^{\ast}(h^{-1}\circ\eta\circ h^{\ast}(\partial
_{a})\lrcorner_{g}R(a\wedge b))\nonumber\\
&  =\underline{\eta\circ h}^{\ast}(g^{-1}(\partial_{a})\lrcorner_{g}R(a\wedge
b))\nonumber\\
&  =\underline{\eta\circ h}^{\ast}\circ\text{ }_{g}R(b),\nonumber\\
h^{\ast}(\partial_{a})\lrcorner\mathcal{R}(a\wedge b)  &  =\mathcal{R}(b).
\label{GEF.11}%
\end{align}

Then, it is obvious that
\begin{equation}
h^{*}(\partial_{b})\cdot\mathcal{R}(b)=\underline{h}^{*}(\partial_{a}%
\wedge\partial_{b})\cdot\mathcal{R}(a\wedge b). \label{GEF.12}%
\end{equation}

Finally, using once again the \emph{master formula} $g^{-1}=h^{-1}\circ
\eta\circ h^{\ast},$ and Eq.(\ref{RRF.6c}), we get
\begin{align}
\underline{h}^{\ast}(\partial_{a}\wedge\partial_{b})\cdot\mathcal{R}(a\wedge
b)  &  =\underline{h}^{\ast}(\partial_{a}\wedge\partial_{b})\cdot
\underline{\eta\circ h}^{\ast}\circ\text{ }_{g}R(a\wedge b)\nonumber\\
&  =\underline{h^{-1}\circ\eta\circ h}^{\ast}(\partial_{a}\wedge\partial
_{b})\cdot\text{ }_{g}R(a\wedge b)\nonumber\\
&  =\underline{g}^{-1}(\partial_{a}\wedge\partial_{b})\cdot\text{ }%
_{g}R(a\wedge b),\nonumber\\
\underline{h}^{\ast}(\partial_{a}\wedge\partial_{b})\cdot\mathcal{R}(a\wedge
b)  &  =\text{ }_{g}R. \label{GEF.13}%
\end{align}

The above results implies that the \emph{Ricci curvature scalar field} for
$(U,\gamma,g)$ can be expressed by the noticeable formula
\begin{align}
_{g}R  &  =\underline{h}^{\ast}(\partial_{a}\wedge\partial_{b})\cdot
((a\cdot\partial_{o}\Omega)(b)-(b\cdot\partial_{o}\Omega)(a)+\Omega
(a)\underset{\eta}{\times}\Omega(b))\label{GEF.13a}\\
& \nonumber
\end{align}

Now, concerning the torsion of the metrical compatible connections $_{\eta
}\mathcal{D}$ and $_{g}\mathcal{D}$, we have putting%
\begin{align}
T^{^{\prime}}(h^{-1}a,h^{-1}b)  &  =\text{ }_{g}\mathcal{D}_{h^{-1}a}%
^{+}h^{-1}b-\text{ }_{g}\mathcal{D}_{h^{-1}b}^{+}h^{-1}a-[h^{-1}%
a,h^{-1}b],\nonumber\\
T(a,b)  &  =\text{ }_{\eta}\mathcal{D}_{a}^{+}b-\text{ }_{\eta}\mathcal{D}%
_{b}^{+}a-[a,b], \label{GEF14}%
\end{align}
that
\begin{equation}
T^{^{\prime}}(h^{-1}a,h^{-1}b)=h^{-1}\left(  T(a,b)\right)  +h^{-1}%
([a,b])-[h^{-1}a,h^{-1}b], \label{GEF15}%
\end{equation}
which shows the important result that even if the connection $_{\eta
}\mathcal{D}$ has null torsion, in general the deformed connection
$_{g}\mathcal{D}$ \ has a \textit{non} null torsion.

\subsection{Deformations Induced by Diffeomorphisms}

We now investigate the what happens when the deformation tensor is induced by
a diffeomorphism. To be precise, consider the structures $(U,\gamma^{\prime
},\eta)$ and $(U,\gamma,g)$. Consider a diffeomorphism
\begin{align}
\mathtt{h} &  :M\rightarrow M,\text{ }\nonumber\\
\mathfrak{e}^{\prime} &  =\mathtt{h}\mathfrak{e}\label{diff1}%
\end{align}
such that the metric tensors $\mathtt{\eta}$ and $%
\slg
$ associated to the metric extensors $\eta$ and $g$ are related by
\begin{equation}%
\slg
=\mathtt{h}^{\ast}%
\mbox{\boldmath{$\eta$}}%
.\label{diff2}%
\end{equation}
The respective $%
\slg
$ and $\mathtt{\eta}$ compatible covariant derivatives are related by the
usual definition\footnote{Or equivalently $\mathtt{h}^{\ast}(_{\mathtt{\eta}%
}\mathcal{D}_{v}u)=$ $_{%
\slg
}\mathcal{D}_{\mathtt{h}^{\ast}v}\mathtt{h}^{\ast}u$ , where $\mathtt{h}%
^{\ast}$ is the pullback operator ($\mathtt{h}^{\ast}u=\mathtt{h}_{\ast}%
^{-1}u$) and now $_{%
\slg
}\mathcal{D}=\mathtt{h}^{\ast}$ $_{\mathtt{\eta}}\mathcal{D}$ is the pullback
of $_{\mathtt{\eta}}\mathcal{D}$ (see, e.g., \cite{rodoliv2006}).}%
\begin{equation}
\left.  (_{\mathtt{\eta}}\mathcal{D}_{v}u)\right\vert _{\mathtt{h}%
\mathfrak{e}}f=\left.  \text{ }_{%
\slg
}\mathcal{D}_{\mathtt{h}_{\ast}^{-1}v}\mathtt{h}_{\ast}^{-1}u\right\vert
_{\mathfrak{e}}f\circ\mathtt{h},\label{diff3}%
\end{equation}
where $_{%
\slg
}\mathcal{D}$ \ is the inverse pushforward of ${}_{%
\slg
}\mathcal{D}$ :$=\mathtt{h}_{\ast\mathtt{\ \eta}}^{-1}\mathcal{D=}%
\mathtt{h}^{\ast}\mathcal{D}$. From Eq.(\ref{diff3}) we have
\begin{equation}
\mathtt{h}_{\ast}^{-1}(_{\mathtt{\eta}}\mathcal{D}_{v}u)=\text{ }_{%
\slg
}\mathcal{D}_{\mathtt{h}_{\ast}^{-1}v}\mathtt{h}_{\ast}^{-1}u\text{
,}\label{df3}%
\end{equation}
since by definition $\left.  (_{\mathtt{\eta}}\mathcal{D}_{v}u)\right\vert
_{\mathtt{h}\mathfrak{e}}f=$ $\left.  \mathtt{h}_{\ast}^{-1}(_{\mathtt{\eta}%
}\mathcal{D}_{v}u)\right\vert _{\mathfrak{e}}f\circ\mathtt{h}$

Let $U$, $\mathtt{h}(U)\subset U_{o}\subset M$ and $\{%
\slx
^{\mu}\}$ be the coordinate functions of \ a local chart $(U_{o},\phi_{o})$ of
a given atlas of $M.$

Let moreover,
\begin{equation}%
\slx
^{\mu}(\mathfrak{e})=x^{\mu}\text{, }%
\slx
^{\mu}(\mathtt{h}\mathfrak{e})=x^{\prime\mu} \label{diff3a}%
\end{equation}
and
\begin{equation}
x^{\prime\mu}=\mathtt{h}^{\mu}(x^{\alpha}), \label{diff3b}%
\end{equation}
where the $\mathtt{h}^{\mu}$ are invertible functions. In the canonical space
$\mathcal{U}_{0}$ (\cite{4} \ the representative of the diffeomorphism
$\mathtt{h}$ is the (vector valued) \textit{function}%
\begin{equation}
\mathfrak{h}:\mathcal{U}\rightarrow\mathcal{U},\text{ }x\mapsto x^{\prime
}=\mathfrak{h}(x), \label{rep diff}%
\end{equation}
\ where $x$ and $x^{\prime}$ are the position vectors of the events
$\mathfrak{e}$ and $\mathfrak{e}^{\prime}$. Let $\mathbf{a}=a^{\alpha}(%
\slx
^{\mu}(\mathfrak{e}))\frac{\partial}{\partial x^{\alpha}}\in\sec TU$ be a
generic vector field. Its representation in$%
{\displaystyle\bigwedge\nolimits^{1}}
\mathcal{U}$ will be denoted by $a$. Also, if \texttt{h }is a diffeomorphism
we represent in$%
{\displaystyle\bigwedge\nolimits^{1}}
$ $\mathcal{U}$ by $\mathtt{h}_{\ast}^{-1}a$ the vector field $\mathtt{h}%
_{\ast}^{-1}\mathbf{a}=a^{\alpha}(%
\slx
^{\mu}(\mathtt{h}\mathfrak{e}))\frac{\partial x^{\mu}}{\partial x^{\prime
\alpha}}\frac{\partial}{\partial x^{\mu}}\in\sec TU$ . Now, associated to the
diffeomorphism \texttt{h} and its representative $\mathfrak{h}$ in the
canonical space $\mathcal{U}_{o}$ there is a unique extensor field $h:%
{\displaystyle\bigwedge\nolimits^{1}}
\mathcal{U}_{0}\rightarrow%
{\displaystyle\bigwedge\nolimits^{1}}
\mathcal{U}_{0}$ such that
\begin{equation}
h(\mathtt{h}_{\ast}^{-1}a)=a. \label{diff4a}%
\end{equation}
Then it is trivial to verify that we can write the extensor field $h$ in terms
of the representatives $\mathfrak{h}$ of \texttt{h }and $\mathtt{h}_{\ast
}^{-1}a$ of $\mathtt{h}_{\ast}^{-1}a$ as
\begin{equation}
h(\mathtt{h}_{\ast}^{-1}a)=(\mathtt{h}_{\ast}^{-1}a)\cdot\mathbf{\partial}%
_{o}\mathfrak{h}(x) \label{diff4}%
\end{equation}
where $(\mathtt{h}_{\ast}^{-1}a)\cdot\mathbf{\partial}_{o}$ is the standard
directional derivative in the direction of the vector field $\mathtt{h}_{\ast
}^{-1}a$. Also, $\partial_{o}\equiv\partial_{ox}$ is the vector derivative
operator in $%
{\displaystyle\bigwedge\nolimits^{1}}
\mathcal{U}$ at $x$ $\in\mathcal{U}$.

Eq.(\ref{diff4a}) and Eq.(\ref{diff4}) shows that the representative of the
field $\mathtt{h}_{\ast}^{-1}\mathbf{a}$ $\in\sec TU$ in $%
{\displaystyle\bigwedge\nolimits^{1}}
\mathcal{U}$ can be written also as $h^{-1}(a)$ which we write simply as
$h^{-1}a$.

On the other hand if
\slg
$=\mathtt{h}^{\ast}%
\mbox{\boldmath{$\eta$}}%
$ we have immediately that the corresponding extensor fields may be related
by\footnote{Of course, if $h^{\prime}=h\Lambda$, where $\Lambda$ is a Lorentz
extensor ($\Lambda\eta\Lambda^{\dagger}=\Lambda^{\dagger}\eta\Lambda=\eta$) it
determines the same $g$ as $h$. We fix a gauge by choosing $\Lambda=1$.}
\begin{equation}
h\circ\eta\circ h^{\dagger}=g,\label{diff6}%
\end{equation}
where\ $h=h^{\dagger}$.

We now prove that for any $a,b\in%
{\displaystyle\bigwedge\nolimits^{1}}
\mathfrak{h}(\mathcal{U)}$ we have that the extensor field $h$ associated to a
diffeomorphism \texttt{h} and its pullback \texttt{h}$^{\ast}$ satisfy:%
\begin{equation}
h[a,b]=[h(a),h(b)], \label{diff7}%
\end{equation}
where $[$ $,$ $]$ is the commutator of vector fields as defined in Section
3.2. of \cite{4}. To start, recall that we must have that
\begin{equation}
\partial_{o}=h(\partial_{o}^{\prime}), \label{diff8}%
\end{equation}
with $\partial_{o}^{\prime}\equiv\partial_{ox^{\prime}}$ the vector operator
at $x^{\prime}$.

Indeed, take any scalar function $f:$ $\mathcal{U}\rightarrow\mathbb{R}$. Let
us calculate
\begin{equation}
a\cdot\partial_{ox}f(\mathfrak{h}(x)). \label{diff81}%
\end{equation}

Using the definition of directional derivative, we have%
\begin{align}
a\cdot\partial_{ox}f(\mathfrak{h}(x))  &  =\lim_{\lambda\rightarrow0}%
\frac{f(\mathfrak{h}(x+\lambda a))-f(\mathfrak{h}(x))}{\lambda}\nonumber\\
&  =\lim_{\lambda\rightarrow0}\frac{f(\mathfrak{h}(x)+\lambda
h(a))-f(\mathfrak{h}(x))}{\lambda}\nonumber\\
&  =h(a)\cdot\partial_{ox^{\prime}}f(\mathfrak{h}(x))=h(a)\cdot\partial
_{ox^{\prime}}f(x^{\prime}). \label{diff82}%
\end{align}
Using moreover the fact that $h=h^{\dagger}$ we have $h(a)\cdot\partial
_{ox^{\prime}}f(x^{\prime})=a\cdot h(\partial_{ox^{\prime}})f(\mathfrak{h}%
(x))=a\cdot\partial_{ox}f(x^{\prime})$ and Eq.(\ref{diff8}) is proved.

Now, using the definition of the commutator of vector fields at $%
{\displaystyle\bigwedge\nolimits^{1}}
\mathfrak{h}(\mathcal{U})$ we have%
\begin{align}
\left.  \lbrack h(a),h(b)]\right\vert _{x^{\prime}}  &  =h(a)\cdot\partial
_{o}^{\prime}h(b)-h(b)\cdot\partial_{o}^{\prime}h(a)\nonumber\\
&  =a\cdot h(\partial_{o}^{\prime})h(b)-b\cdot h(\partial_{o}^{\prime
})h(a)\nonumber\\
&  =a\cdot\partial_{o}h(b)-b\cdot\partial_{o}h(a)\nonumber\\
&  =h(a\cdot\partial_{o}b-b\cdot\partial_{o}a)\nonumber\\
&  =\left.  h[a,b]\right\vert _{x}. \label{diff9}%
\end{align}
where we use Eq.(\ref{diff8}) and the properties of the adjoint.

With this result we see that if $h$ is generated by a diffeomorphism as above,
then $h^{-1}([a,b])-[h^{-1}a,h^{-1}b]=0$ and we conclude the important result
that for diffeomorphism \texttt{h}, if the structure $(M,\gamma^{\prime}%
,\eta)$ has zero curvature and torsion tensors, also the deformed structure
$(M,\gamma,g)$ has zero curvature and torsion tensors.

\section{Levi-Civita Geometric Structure}

Let $(U,\lambda,g)$\emph{\ }be the \emph{Levi-Civita geometric structure} on
$U$.\emph{\ }We recall \cite{5} that $\lambda$ is the \emph{Levi-Civita
connection field} which is given by
\begin{equation}
\lambda(a,b)=\frac{1}{2}g^{-1}\circ(a\cdot\partial_{o}g)(b)+\omega
_{0}(a)\underset{g}{\times}b, \label{LGS.0a}%
\end{equation}
where $\omega_{0}$ is the smooth $(1,2)$-extensor field given by
\begin{equation}
\omega_{0}(a)=-\frac{1}{4}\underline{g}^{-1}(\partial_{b}\wedge\partial
_{c})a\cdot((b\cdot\partial_{o}g)(c)-(c\cdot\partial_{o}g)(b)). \label{LGS.0b}%
\end{equation}

As we know, such $\lambda$ is \emph{symmetric} by interchanging of their
vector variables, i.e., $\lambda(a,b)=\lambda(b,a).$ The $g$-\emph{compatible}
pair of $a$-\emph{DCDO's }associated to $(U,\lambda,g)$\emph{\ }is, of course,
the so-called pair of \emph{Levi-Civita }$a$-\emph{DCDO's,} namely $(D_{a}%
^{+},D_{a}^{-}).$

We now present the basic properties which are satisfied by the curvature field
of $(U,\lambda,g).$\vspace{0.1in}

\textbf{i.} The curvature field $\rho$ has the \emph{cyclic property}
\begin{equation}
\rho(a,b,c)+\rho(b,c,a)+\rho(c,a,b)=0. \label{LGS.1a}%
\end{equation}

\textbf{ii.} The curvature field $\rho$ satisfies the so-called
\emph{Bianchi's identity}, i.e.,
\begin{equation}
\nabla_{d}^{+++-}\rho(a,b,c)+\nabla_{a}^{+++-}\rho(b,d,c)+\nabla_{b}%
^{+++-}\rho(d,a,c)=0. \label{LGS.1aa}%
\end{equation}

We emphasize that these remarkable properties hold for any symmetric
parallelism structure, a parallelism structure $(U,\gamma)$ in which the
connection field $\gamma$ is\emph{\ symmetric}, i.e., $\gamma(a,b)=\gamma
(b,a),$ see \cite{5}.

\textbf{iii.} The curvature field $\rho$ satisfies also the symmetry property
\begin{equation}
\rho(a,b,c)\underset{g}{\cdot}d=\rho(c,d,a)\underset{g}{\cdot}b.
\label{LGS.1b}%
\end{equation}

To prove Eq.(\ref{LGS.1b}) we use four times Eq.(\ref{LGS.1a}) obtaining
\begin{align*}
\rho(a,b,c)\underset{g}{\cdot}d+\rho(b,c,a)\underset{g}{\cdot}d+\rho
(c,a,b)\underset{g}{\cdot}d  &  =0,\\
\rho(b,c,d)\underset{g}{\cdot}a+\rho(c,d,b)\underset{g}{\cdot}a+\rho
(d,b,c)\underset{g}{\cdot}a  &  =0,\\
-\rho(c,d,a)\underset{g}{\cdot}b-\rho(d,a,c)\underset{g}{\cdot}b-\rho
(a,c,d)\underset{g}{\cdot}b  &  =0,\\
-\rho(d,a,b)\underset{g}{\cdot}c-\rho(a,b,d)\underset{g}{\cdot}c-\rho
(b,d,a)\underset{g}{\cdot}c  &  =0.
\end{align*}
Now, by adding the above equations , and using six times Eq.(\ref{RRF.2b}) and
twice Eq.(\ref{RRF.2}), we get
\[
2\rho(a,b,c)\underset{g}{\cdot}d-2\rho(c,d,a)\underset{g}{\cdot}b=0.
\]
Whence, the expected result immediately follows.

\textbf{iv.} The \emph{Riemann field} for $(U,\lambda,g)$,\emph{\ }namely
$B\mapsto R(B),$ satisfies the following cyclic property
\begin{equation}
R(a\wedge b)\times c+R(b\wedge c)\times a+R(c\wedge a)\times b=0.
\label{LGS.1c}%
\end{equation}

To prove Eq.(\ref{LGS.1c}) we use three times Eq.(\ref{RRF.6a}), obtaining
\begin{align*}
R(a\wedge b)\times c  &  =g\circ\rho(a,b,c),\\
R(b\wedge c)\times a  &  =g\circ\rho(b,c,a),\\
R(c\wedge a)\times b  &  =g\circ\rho(c,a,b).
\end{align*}
Then, by adding the above equations, and using Eq.(\ref{LGS.1a}), the result follows.

We finally present the symmetry properties of the \emph{Riemann} and
\emph{Ricci} \emph{fields} corresponding to the Levi-Civita geometric
structure.\vspace{0.1in}

\textbf{i.} The Riemann field is \emph{symmetric}, i.e.,
\begin{equation}
R^{\dagger}(B)=R(B). \label{LGS.2}%
\end{equation}

To prove this result let us take $a,b,c,d\in\mathcal{V}(U).$ Then, using the
fundamental property for the adjoint operator $\left.  {}\right.  ^{\dagger},$
and Eq.(\ref{RRF.3}) and Eq.(\ref{LGS.1b}), we can write
\begin{align*}
R^{\dagger}(a\wedge b)\cdot(c\wedge d)  &  =(a\wedge b)\cdot R(c\wedge d)\\
&  =-\rho(c,d,a)\underset{g}{\cdot}b=-\rho(a,b,c)\underset{g}{\cdot}d,\\
&  =R(a\wedge b)\cdot(c\wedge d).
\end{align*}

Now, by applying the \emph{bivector operator} $\dfrac{1}{2}\partial_{c}%
\wedge\partial_{d}$ on both sides of equation above, we have
\begin{align*}
\dfrac{1}{2}\partial_{c}\wedge\partial_{d}R^{\dagger}(a\wedge b)\cdot(c\wedge
d)  &  =\dfrac{1}{2}\partial_{c}\wedge\partial_{d}R(a\wedge b)\cdot(c\wedge
d),\\
R^{\dagger}(a\wedge b)  &  =R(a\wedge b).
\end{align*}

But, by applying the \emph{scalar operator} $\dfrac{1}{2}B\cdot(\partial
_{a}\wedge\partial_{b})$ on both sides of this equation, we get
\begin{align*}
\dfrac{1}{2}B\cdot(\partial_{a}\wedge\partial_{b})R^{\dagger}(a\wedge b)  &
=\dfrac{1}{2}B\cdot(\partial_{a}\wedge\partial_{b})R(a\wedge b)\\
R^{\dagger}(\dfrac{1}{2}B\cdot(\partial_{a}\wedge\partial_{b})a\wedge b)  &
=R(\dfrac{1}{2}B\cdot(\partial_{a}\wedge\partial_{b})a\wedge b),\\
R^{\dagger}(B)  &  =R(B).
\end{align*}

\textbf{ii.} The Ricci field is \emph{symmetric,} i.e.,
\begin{equation}
R^{\dagger}(b)=R(b). \label{LGS.3}%
\end{equation}

To prove this result let us take $b,c\in\mathcal{V}(U).$ A straightforward
calculation using Eq.(\ref{RRF.5}) and Eq.(\ref{LGS.2}) allows us to get
\begin{align*}
R^{\dagger}(b)\cdot c  &  =b\cdot R(c)\\
&  =b\cdot(g^{-1}(\partial_{a})\lrcorner R(a\wedge c))=(g^{-1}(\partial
_{a})\wedge b)\cdot R(a\wedge c)\\
&  =R(g^{-1}(\partial_{a})\wedge c)\cdot(a\wedge b)=(g^{-1}(\partial
_{a})\wedge c)\cdot R^{\dagger}(a\wedge b)\\
&  =c\cdot(g^{-1}(\partial_{a})\lrcorner R^{\dagger}(a\wedge b))=c\cdot
(g^{-1}(\partial_{a})\lrcorner R(a\wedge b))\\
&  =c\cdot R(b),
\end{align*}
hence, by the non-degeneracy of the scalar product, the required result follows.

\section{Conclusions}

Here, the last paper in a series of four where we end our presentation of the
basics of a systematical approach of the geometric algebra of multivector and
extensors to the differential geometry of an arbitrary smooth manifold $M$
equipped with a metric tensor
\slg
and a connection $\nabla$. The theory of the Riemann and Ricci fields
associated to the triple $(M,\nabla,%
\slg
)$ is studied in an arbitrary open set $U\subset M$ through the introduction
of a geometric structure on $U$, i.e., a triple $(U,\gamma,g)$ where $\gamma$
is a general connection field on $U$ (describing there the "effects" of
$\nabla$) and $g$ is a metric extensor field associated to $%
\slg
$. The relation between geometrical structures related by gauge extensor
fields has been clarified. These geometries may be said to be deformations one
of each other. Moreover, we studied in details the important case of a class
of deformed Levi-Civita geometrical structures and proved some key theorems
about these structures that are important in the formulation of geometric
theories of the gravitational field. Our main result is that if we start with
a \textit{flat} connection in $M$, represented in $U$ by $(U,\gamma^{\prime
},\eta)$ the deformed geometry structure $(U,\gamma,g)$ in general can only be
the representiative on $U$ of a connection on $M$ which has \textit{null}
curvature but a \textit{non null }torsion. Only in the case in which the
deformation tensor $h$ is associated to a diffeomorphism in a precise sense
defined in section 2.3 is that starting with $(U,\gamma^{\prime},\eta)$, with
null curvature and torsion tensors the deformed geometry $(U,\gamma,g)$ also
has null curvature and torsion tensors. Of course, there are additional
related topics that merit a presentation with our formalism and this will be
done opportunely.\medskip

\textbf{Acknowledgments:} V. V. Fern\'{a}ndez and A. M. Moya are very grateful
to Mrs. Rosa I. Fern\'{a}ndez who gave to them material and spiritual support
at the starting time of their research work. This paper could not have been
written without her inestimable help. Authors are also grateful to Drs. E.
Notte-Cuello and E. Capelas de Oliveira for useful discussions.

\end{document}